\documentclass{article}
\usepackage{amsmath,amssymb,amsthm}
\title{Reals $n$-generic relative to some perfect tree}
\author{Bernard A. Anderson\\Department of Mathematics\\University of California\\Berkeley, CA 94720, USA}
\newtheorem{theorem}{Theorem}[section]
\newtheorem{lemma}[theorem]{Lemma}
\newtheorem{corr}[theorem]{Corollary}
\newtheorem{prop}[theorem]{Proposition}
\newtheorem*{claim}{Claim}
\newcommand{\KlO}{\ensuremath{\mathcal{O}}}
\newcommand{\bor}{\ensuremath{\mathcal{B}}}
\newcommand{\WF}{\ensuremath{W \! F}}
\newcommand{\setsep}{\ensuremath{\, | \;}}
\begin{document}
\maketitle
\begin{abstract}
We say that a real $X$ is $n$-generic relative to a perfect tree $T$ if $X$ is a path through $T$ and for all 
$\Sigma^0_n (T)$ sets $S$, there exists a number $k$ such that either $X|k \in S$ or for all $\sigma \in T$
extending $X|k$ we have $\sigma \notin S$\@.  A real $X$ is $n$-generic relative to some perfect tree if there exists
such a $T$.  We first show that for every number $n$ all but countably many reals are $n$-generic relative to some perfect
tree.  Second, 
we show that proving this statement requires ZFC$^-$ + ``$\exists$ infinitely many iterates of the power set of $\omega$''.  
Third, we prove that every finite iterate of the hyperjump, $\KlO^{(n)}$, is not 2-generic relative to any perfect tree and for every 
ordinal $\alpha$ below the least $\lambda$ such that $\sup_{\beta < \lambda} (\beta$th admissible$) = \lambda$, the iterated hyperjump 
$\KlO^{(\alpha)}$ is not 5-generic relative to any perfect tree.  Finally, we demonstrate some necessary conditions for 
reals to be 1-generic relative to some perfect tree.
\end{abstract}
\section{Introduction}
A real (viewed as an element of $2^\omega$) is $n$-generic if for every $\Sigma^0_n$ set there is an initial segment of the real 
which either meets the set or for which no extension of the segment can meet the set.  These reals have many interesting 
char\-acter\-istics and have been studied extensively (see Jockusch and Posner \cite{generic2} and Kumabe \cite{generic1} among 
others).  
While the set of $n$-generics is comeager, it is in some ways limited.  In particular, it is completely excluded from the 
cone above $\mathbf{0^\prime}$ since no 1-generic can compute an r.e.\ set. 

A question which naturally arises then is how this set might be expanded from reals which are $n$-generic to those that can be 
made to seem $n$-generic in some appropriate context.  An attractive framework for this question is to consider reals which are 
$n$-generic when viewed as paths through a given perfect tree, rather than all of $2^\omega$.

\medskip\noindent\textbf{Definition:} A real $X$ is $n$-generic relative to a perfect tree $T$ if $X$ is a path through $T$ and 
for all $\Sigma^0_n (T)$ sets $S$, there is a $k$ such that either $X|k \in S$ or $\sigma \notin S$ for every  $\sigma \in T$ 
extending $X|k$.

\medskip\noindent\textbf{Definition:} A real $X$ is $n$-generic relative to some perfect tree if there exists a perfect tree $T$ 
such that $X$ is $n$-generic relative to $T$.

\medskip
This definition results in a version of genericity which includes many reals which seem ``essentially'' generic.  For example, if $G$ 
is a $n$-generic real and $R$ is recursive then neither $G \oplus R$ nor $G \oplus G$ is $n$-generic, but both are $n$-generic relative to 
some perfect tree.  Another form of relative genericity that has been studied is genericity relative to another real.  Given a real $A$ 
we say that a real $X$ is $n$-generic($A$) if for every $\Sigma^0_n (A)$ set $S$, either $X$ meets $S$ or there is an initial segment of $X$ 
which has no extensions in $S$.  This results in a more restrictive version of genericity which is similar to higher levels of genericity.  
For example a real is 1-generic($0^\prime$) if and only if it is 2-generic.  We note that every real which is $n$-generic($A$) for some 
real $A$ is $n$-generic and every real which is $n$-generic is $n$-generic relative to some perfect tree.

In this paper, we examine the set of reals $n$-generic relative to some perfect tree.  We first show that the set of reals not $n$-generic 
relative to any perfect tree is countable.  From this we can infer that many reals with properties that are not normally associated with 
genericity still seem generic in the context of some perfect tree.  For example, there are reals of minimal degree and reals with high 
information content, such as the theory of second order arithmetic, that are generic relative to some perfect tree.  

The proof that the set of reals not $n$-generic relative to any perfect tree is countable uses ZFC$^-$ and $n$ iterates of the power 
set of $\omega$.  We show that for sufficiently large $n$, this requirement is sharp and cannot be significantly improved.  From this we 
see that for reasonably high values of $n$, the set of reals not $n$-generic relative to any perfect tree is unusually large (rich) 
for a countable set of this type.  It provides a natural example of a set which needs this level of ZFC to be understood.
  
While the result above holds only for sufficiently high values of $n$, for lower values the set of reals not $n$-generic relative to 
any perfect tree is also rich.  By looking at the iterates of the hyperjump, we demonstrate that the set still contains reals of 
unexpectedly high complexity.  Even for $n=2$, relatively large fragments of arithmetic fail to prove the set is countable.  We also 
begin to characterize the sets that are 1-generic relative to some perfect tree.

These results are in a similar vein to results of Reimann and Slaman \cite{Random}, who have studied the set of reals which appear 
random in some context, in this
case relative to some continuous measure.  Our results for genericity are analogous to what they 
discovered for randomness in surprisingly many, but not all, instances.

This work formed part of the author's Ph.D. Thesis at the University of California at Berkeley.  We thank Theodore Slaman, the dissertation 
supervisor, for his introduction to the topic and his repeated suggestions of new approaches to problems.
\section{Co-Countably Many Reals}
We wish to show that the set of reals not $n$-generic relative to any perfect tree is countable.  D. Martin \cite{Martin} used 
Borel determinacy to show that any property which is Borel and cofinal in the Turing degrees is represented on every degree in a 
cone of Turing degrees.  The base of this cone is the complexity of the winning strategy for an associated game.  

\begin{theorem}[Martin \cite{Martin}] Let $\bor$ be a Borel set of reals such that for every Turing degree $\mathbf{d}$ there is an 
$\mathbf{e} \geq_T \mathbf{d}$ and an $X$ in $\mathbf{e}$ such that $X \in \bor$.  Then there is a degree $\mathbf{c}$ such that 
for all $\mathbf{b} \geq_T \mathbf{c}$ there is a $Y$ in $\mathbf{b}$ such that $Y \in \bor$. \label{Marlem} \end{theorem}
\begin{proof} Consider a two person game where player I constructs a real $X$ and player II constructs a real $Y$.  Play alternates 
between the players, each adding the next digit to the real they are constructing for their turn.  Player I wins iff $Y \leq_T X$ and 
either $X \not\le_T Y$ or $X \in \bor$.  

By Borel determinacy, there exists a winning strategy $\sigma$.  Suppose $\sigma$ is a winning strategy for II.  Since $\bor$ is 
cofinal in the Turing degrees, let $Z \in \bor$ with $\sigma \leq_T Z$\@.  Let I play $Z$ and II play according to $\sigma$ resulting 
in 
$Y$\@.  But then since $\sigma \leq_T Z$, we have $Y \leq_T Z$ with $Z \in \bor$ so I wins for a contradiction.  Hence $\sigma$ 
must be 
a winning strategy for I.
 
Let $Z \geq_T \sigma$ be an arbitrary real in the cone above $\sigma$.  Have II play $Z$ and I play according to $\sigma$, 
resulting in $X$\@.  Since $\sigma \leq_T Z$ we have $X \leq_T Z$\@.  Since I wins, $Z \leq_T X$ and $X \in \bor$.  Hence for any $Z$ 
in 
the cone, there is an $X \in \bor$ such that $X \equiv_T Z$\@.  
\end{proof} 
Reimann and Slaman \cite{Random} have developed a powerful way to relativize this lemma.  Let $\bor \subseteq 2^\omega \times 2^\omega$
denote a set of reals where the first real holds some property relative to the second.  Let $\bor^Z = \{X \setsep (X,Z) \in \bor \}$ and let 
the 
notation $X \equiv_{T,A} Y$ mean $X \oplus A \equiv_T Y \oplus A$\@.  Suppose that for every $Z$ the set $\bor^Z$ is Borel in $Z$ 
and cofinal in the Turing degrees as in the above method for generating a cone.  They prove that for all but countably many reals $X$, 
there exist reals $Y$ and $G$ such that $X \equiv_{T,G} Y$ and $Y \in \bor^G$\@.  We outline this proof in the next paragraph.  

Let $\beta$ be the least ordinal such that $L_\beta$ satisfies enough ZFC ($L_\beta$ is countable) and let \mbox{$X \notin L_\beta$} be 
arbitrary.  Reimann and Slaman use Kumabe-Slaman forcing to find a real $G$ such that $L_\beta [G] \models$ ZFC and every element of 
$2^\omega \cap L_\beta [G]$ is recursive in $X \oplus G$\@.  In particular, the strategy for the game in the proof of Theorem 
\ref{Marlem} relative to $G$ is recursive in $X \oplus G$\@.  So by Theorem \ref{Marlem} relativized to $G$ there exists 
$Y \in \bor^G$ with $Y \equiv_{T,G} X$\@.  

Thus to prove all but countably many reals are $n$-generic relative to some perfect tree, we need to find a set $\bor$ such that for any 
$X, Y, G$ with $Y \equiv_{T,G} X$ and $Y \in \bor^G$ we have $X$ $n$-generic relative to some perfect tree.  $\bor$ must also be Borel and 
such that for every $Z$ the set $\bor^Z$ is cofinal in the Turing degrees.  We find it suffices to let $\bor$ be the set of reals of 
Turing degree $X \oplus A$ for any $X,A$ such that $X$ is ($n+1$)-generic($A$).  We use the following lemma.  

\begin{lemma} Let $n \geq 2$, $A$ be a set, $X$ be $n$-generic($A$), and $X \equiv_{T,A} Y$\@. Then $Y$ is ($n-1$)-generic relative 
to some perfect tree. \label{mainlem}\end{lemma}
\begin{proof} Let $\Psi : X \to Y$ and $\Phi : Y \to X$ be $A$-recursive Turing reductions that witness $X \equiv_{T,A} Y$\@.  Since 
$X$ is at least 2-generic($A$), let $p \in X$ be such that $p \Vdash \Phi \circ \Psi = \mbox{id}\ \wedge\ \Psi$ total (i.e.\ this statement 
holds for all 2-generic($A$) reals extending $p$).  Let 
$T=\{\sigma \setsep \exists q \supseteq p[\sigma \subseteq \Psi (q)] \}$.  $T$ is a perfect tree by our choice of $p$.  We claim that $Y$ is 
($n-1$)-generic relative to $T$\@.  

Let $S$ be an arbitrary $\Sigma^0_{n-1} (T)$ set.  We consider the pullback $\Psi^{-1} (S) = \{x \setsep \exists y[\Psi (x) \supseteq y\ \wedge\ y 
\in S] \}$.  $T$ is $\Sigma^0_1(A)$ so $S$ is $\Sigma^0_n (A)$ and $\Psi^{-1} (S)$ is $\Sigma^0_n (A)$.  We will now apply the 
genericity of $X$ for the pullback to get the genericity of $Y$ for $S$\@.  

Since $X$ is $n$-generic ($A$) we have two possible cases.

\noindent Case 1: $\exists n [X|n \in \Psi^{-1} (S)]$.  We then let $m$ be such that $Y|m \subseteq \Psi (X|n)$ and $Y|m \in S$\@.  

\noindent Case 2: $\exists n \forall q \supseteq X|n [q \notin \Psi^{-1} (S)]$.  Let $m$ be such that  $\Phi (Y|m) \supseteq X|n$.  
We will show $Y|m$ witnesses $Y$ is ($n-1$)-generic relative to $T$ for $S$ (no extension in $T$ of $Y|m$ is in $S$).  Consider 
an arbitrary $r \in T$ such that $r \supseteq Y|m$.  Since $r \in T$, let $q$ be such that $\Psi (q) \supseteq r$ and  $q \supseteq p$.  
We note $q \supseteq \Phi (\Psi (q)) \supseteq \Phi (r) \supseteq \Phi (Y|m) \supseteq X|n$.  Hence by the condition for this case, 
$q \notin \Psi^{-1} (S)$ so $r \notin S$.  Since $r$ is arbitrary, for all $r \supseteq Y|m$ with $r \in T$ we have $r \notin S$\@.  
\end{proof}
We note that a similar proof can be used to show for $n \geq 1$ that sets in the same truth table degree as an $n$-generic are 
$n$-generic relative to some perfect tree.

We can now use the approach outlined above.

\begin{theorem} For every $n\in\omega$, the set of reals not $n$-generic relative to some perfect tree is countable. \end{theorem}
\begin{proof} Fix $n \in \omega$ and let 
$$\bor = \{(x,z) \setsep \exists c \exists h [h \mbox{ is } (n+1) \mbox{-generic}(c \oplus z) \mbox{ and } x \equiv_{T,z} c \oplus h \}$$
$\bor$ is arithmetic (since $c,h \leq_T x \oplus z$) so $\bor$ is Borel.  Given any reals $C$ and $Z$, we let $H$ be 
$(n+1)$-generic($C \oplus Z$) and $X = H \oplus C$ to get $X \in \bor^Z$ with $X \geq_T C$\@.  Hence $\bor^Z$ is cofinal in the Turing 
degrees.  By the theorem of Reimann and Slaman \cite{Random} noted above, for all but countably many reals $Y$, there exist $X$ and $G$ 
such that $Y \equiv_{T,G} X$ and $X \in \bor^G$.  

Given such $X$ and $Y$, there exist reals $A$ and $H$ such that $A \oplus H \equiv_{T,G} X$ and $H$ is $(n+1)$-generic($A \oplus G$).  Hence 
$X \equiv_{T, A \oplus G} H$, so $Y \equiv_{T, A \oplus G} H$\@.  By Lemma \ref{mainlem}, $Y$ is $n$-generic relative to some perfect 
tree.  Therefore, all but countably many reals are $n$-generic relative to some perfect tree.    
\end{proof}
\section{ZFC$^-$ and Infinitely Many Iterates of the \\Power Set of $\omega$ Required}
If we examine the proof that the set of reals not $n$-generic relative to any perfect tree is countable, we see that the greatest use of 
the axioms of ZFC comes from the application of Borel determinacy.  The proof uses determinacy of a $\Pi^0_{n+3}$ game on 
$\omega^\omega$, so it requires ZFC$^-$ and the existence of $n$ iterates of the power set of $\omega$ \cite{MarBook}.  We prove that for 
sufficiently large $n$ this is essentially the best possible result.  As a consequence, we show that for any finite $k$ the statement 
``For all $n$, the set of reals not $n$-generic relative to any perfect tree is countable'' cannot be proved from ZFC$^-$ and $k$ 
iterates of the power set of $\omega$.  This suggests the set of reals not $n$-generic relative to any perfect tree is a countable set 
of considerable size and complexity.

\begin{theorem} For every $k \in \omega$ the statement ``For all $n$, the set of reals not $n$-generic relative to any perfect tree is 
countable'' cannot be proved from ZFC$^- + $\mbox{``$ \exists k$ iterates of the power set of $\omega$''.} \end{theorem}

To prove this theorem we use a template developed by Reimann and Slaman \cite{Random} for reals random relative to a continuous 
measure.  We work with the case $k=0$; the general case follows the same pattern.  Let $\lambda$ be the least ordinal such that 
$L_\lambda \models \mbox{ZFC}^-$ and let $O$ be the set of limit ordinals below $\lambda$.  Let $M_\alpha$, for $\alpha \in O$, denote 
master codes.  These are the elementary diagrams of canonical countings of $L_\alpha$.  Reimann and Slaman prove the theorem by showing 
that for some fixed $n$, for every $\alpha \in O$, the master code $M_\alpha$ is not $n$-random relative to a continuous measure.  Since 
the set of master codes is not a countable set in $L_\lambda$, we then have $L_\lambda$ does not satisfy ``For all $n$, the set of reals 
not $n$-random relative to any continuous measure is countable'' but $L_\lambda \models \mbox{ZFC}^-$.  Thus the statement cannot be 
proved from ZFC$^-$.

To show this, they assume towards a contradiction that some $M_\beta$ is $n$-random relative to the measure $\mu$.  It is arithmetic 
to say that $M$ is a master code for an $\omega$-model of ``$V=L_\alpha$ and $\alpha$ a limit and $\alpha \not\geq \lambda$''.  Note such 
an $\omega$-model need not be well-founded.  They show it is also arithmetic to require that there exists a fixed $m \in \omega$ such that 
for all such $M$ and $N$ either one coded model embeds into the other or there is a $\Sigma^0_m (M \oplus N)$ set witnessing the ill-foundedness 
of one of the coded models.  

Reimann and Slaman define a set $\mathcal{M}$, arithmetic in $\mu$, of such psuedo-master codes which are recursive in $\mu$ and not 
shown to be ill-founded by such a comparison.  They then define an order on $\mathcal{M}$ such that the well-founded part of this 
order, $I$, is arithmetic in $\mu \oplus M_\beta$ and equals the set of $M \in \mathcal{M}$ which are actual master codes $M_\alpha$.  
Since random sets cannot accelerate the calculation of well-foundedness, $I$ is arithmetic in $\mu$.

Let $\gamma \leq \beta$ be least such that $M_\gamma \not\leq_T \mu$.  Since $\gamma < \lambda$ there is a real 
$X \in \mbox{Def}(L_\gamma)\setminus L_\gamma$.  By taking a Skolem hull of the parameters defining $X$, Reimann and Slaman show that 
$M_\gamma$ is arithmetic in $I$, hence arithmetic in $\mu$, and $M_\gamma \leq_T M_\beta$.  Since randomness cannot accelerate arithmetic 
definability, $M_\gamma \leq_T \mu$ for a contradiction.

This proof uses only two facts about randomness.  Namely, that it cannot accelerate arithmetic definability or calculations of 
well-foundedness.  We can give an abstract summary as follows.  Let $R(G,T,n)$ be a $\Delta_1^1$ predicate and suppose that for all 
numbers $l$ there exists a number $n$ such that for all $G$ and $T$ such that $R(G,T,n)$ and all $k,m<l$ the following two statements 
hold.  First, for any real $A$, if $A$ is $\Sigma^0_k (T)$ and $\Sigma^0_m (G \oplus T)$ then $A$ is $\Sigma^0_m (T)$.  Second, if $L$ 
is a linear order and $\WF$ the well founded part of $L$ such that $L$ is $\Delta^0_k (T)$ and $\WF$ is $\Delta^0_k (G \oplus T)$ then 
$\WF$ is $\Delta^0_k (T)$.  We can then conclude that for every number $k$ the statement ``For all $n$, the set of reals $G$ such that 
for no $T$ does $R(G,T,n)$ is countable'' cannot be proved from ZFC$^- + \exists k$ iterates of the power set of $\omega$.  We note 
that for some choices of $R$ the statement will not hold at all, in which case this is trivial.

To complete a similar proof for genericity relative to the perfect tree $T$ in place of randomness relative to the measure $\mu$ we will 
demonstrate the corresponding facts for genericity.  We will show that for any fixed $m$ and $k$, and for $n$ sufficiently large relative 
to $m$ and $k$, if $G$ is $n$-generic relative to the perfect tree $T$, then:

\begin{enumerate}
\item{If $A$ is $\Sigma^0_k (T)$ and $\Sigma^0_m (G \oplus T)$ then $A$ is $\Sigma^0_m (T)$.}
\item{If $\WF$ is the well founded part of a linear order recursive in $T$ and \mbox{$\WF \leq_T G \oplus T$} then $\WF \leq_T T$ 
(weaker than above, but see below).}
\end{enumerate}

To show the first fact, we can routinely relativize to a perfect tree the proof that for reals $A,G$ where $G$ is $k$-generic and $A$ is 
$\Sigma^0_k$ and $\Sigma^0_m (G)$, we get that $A$ is $\Sigma^0_m$.  The second fact suffices for our purposes since if $L$ is 
$\Delta^0_k (T)$ and $G$ is $(n+k-1)$-generic relative to $T$ then $G$ is $n$-generic($L$) relative to $T$.  To prove the second fact, 
we use the following lemma.  

\begin{lemma} Let $T$ be a perfect tree and $L$ a linear order of $\omega$ where $L \leq_T T$.  Let $\WF$ be the well founded part 
of $L$\@.  Let $G$ be 2-generic relative to $T$ and such that $\WF \leq_T G \oplus T$\@.  Then $\WF \leq_T T$\@. \end{lemma} 
\begin{proof} By (1.) above, it suffices to show $\WF$ is $\Sigma^0_2 (T)$\@.  For $c \in \omega$ and $D$ an 
initial segment of $L$ we use the notation $c \in D$ to mean there is an ordered pair in $D$ which contains $c$.  Since $L \leq_T T$, if 
$T$ can compute which numbers are in an initial segment then $T$ can compute the initial segment itself.

Let $\Phi$ be a $T$-recursive Turing reduction such that $\Phi (G) = \WF$\@.  For $b \in L$ let $P(b)$ be the set of reals which code initial 
segments of $L$ below $b$.  Let $R(b) \leq_T L$ be the tree defined below such that $P(b)$ is the set of paths through $R(b)$.  
\begin{eqnarray*} R(b) = \{ \sigma \in 2^{<\omega} \setsep \forall m,n < \mbox{length} (\sigma)[(m \in \sigma \rightarrow m \leq_L b )\ \wedge\ \\ 
((m \in \sigma \ \wedge\ n \leq_L m) 
\rightarrow n \in \sigma)]\} \end{eqnarray*}
We define $Q$ as the set of strings in $T$ below which $\Phi$ does not split on $T$\@.
$$Q = \{ \sigma \in T \setsep \neg \exists \tau, \gamma \in T [\tau, \gamma \supseteq \sigma\ \wedge\ \Phi (\tau) \perp \Phi (\gamma)]\}$$
Suppose for some $n$, $G|n \in Q$\@.  Then we can calculate whether $m \in \WF$ by looking for the first 
$\sigma \in T$ such that $\sigma \supseteq G|n$ and $[\Phi (\sigma)](m)\!\!\downarrow$ and taking its value.  Hence $\WF \leq_T T$ and 
we are done.  Thus we may assume for all $n$, $G|n \notin Q$\@.  Since $G$ is 2-generic relative to $T$, there is an $l$ such that for 
all $\tau \supseteq G|l$ with $\tau \in T$ we have $\tau \notin Q$\@.  

We will use the fact that $b \in \WF$ iff $\WF \notin P(b)$.  Let $S = \{ \sigma \setsep \Phi (\sigma)\!\downarrow\ \wedge\ \Phi (\sigma)\notin R(b) \}$.  
We will determine if $b \in \WF$ by checking for the existence of a real which computes an element of $P(b)$ using $\Phi$ and is generic for $S$ 
(some initial segment has no extension in $S$).  Let $\Theta (b)$ be the statement
$$\exists \sigma \in T [\sigma \supseteq G|l \ \wedge\ \forall \tau \in T[\tau \supseteq \sigma \rightarrow \Phi (\tau) \in R(b)]]$$
\begin{claim} $b \notin \WF \Leftrightarrow \Theta (b)$.  \end{claim}
\begin{proof} $(\Longrightarrow)$ $\WF \in P(b)$ since $b \notin \WF$\@.  Hence $\Phi (G) \in P(b)$ so $G$ does not meet $S$\@.  Since $G$ is 
1-generic relative to $T$, there is a $k$ such that for all $\tau \in T$ with $\tau \supseteq G|k$ we have $\tau \notin S$\@.  Thus 
$G|k$ witnesses $\Theta (b)$.  

$(\Longleftarrow)$ Let $\sigma$ witness $\Theta (b)$.  Then for all $\tau \in T$ with $\tau \supseteq \sigma$ we have 
$\Phi (\tau) \in R(b)$.  Also, since $\sigma \supseteq G|l$, for all such $\tau$ we have $\tau \notin Q$ so $\Phi$ splits on $T$ below 
$\tau$.  Using these facts we can construct a perfect subtree of $R(b)$ by applying $\Phi$ to $T$ below $\sigma$.  Hence $P(b)$ is 
uncountable.  Since there are only countably many well founded initial segments of $L$, $b \notin \WF$\@.  
\end{proof}
By the claim $\WF$ is $\Sigma^0_2 (T)$ as desired.  Hence $\WF \leq_T T$\@. 
\end{proof}
\section{Iterated Hyperjumps}
We now look at the set of reals which are $n$-generic relative to some perfect tree for low values of $n$.  We still find that the set
of reals not $n$-generic relative to any perfect tree is a large countable set.  It contains reals of high complexity and its 
countability cannot be proved in large fragments of second order arithmetic.  We show that the finite iterates of the hyperjump, 
$\KlO^{(n)}$, are not 2-generic relative to any perfect tree and the iterates $\KlO^{(\alpha)}$ are not 5-generic relative to any 
perfect tree for any $\alpha$ below the least $\lambda$ such that $\sup_{\beta < \lambda} (\beta$th admissible$) = \lambda$. 

We start with an outline of the proof for the case of $\KlO$. This set can be viewed as $\{e \setsep U_e $ is well-founded$\}$ where 
$U_e$ denotes the $e$th recursive tree in $\omega^{<\omega}$.  We note $\KlO$ then has the property that the well-foundedness of subtrees cannot 
contradict the decision made for the parent tree.  This can be characterized by a $\Sigma^0_2$ set, $S$, so that if $\KlO$ were 
2-generic relative to some $T$ then $T$ would be able to calculate $\KlO$ by tracing subtrees.
 
\begin{lemma} $\KlO$ is not 2-generic relative to any perfect tree. \label{baseO}\end{lemma}
\begin{proof} Suppose not, witnessed by $T$\@.  Recall $\KlO = \{e \setsep U_e $ is well-founded$\}$.  Let $h$ be a recursive 
function defined by $U_{h(e,\gamma)}=\{\sigma \in U_e \setsep \sigma \subseteq \gamma\ \vee\ \sigma \supseteq \gamma \}$.  Let 
\begin{eqnarray*} S=\{\tau \in T \setsep \exists n \exists l [\tau (n)=0 \ \wedge\ \neg (\exists \gamma \in U_n) (\exists \theta \in T) 
[\mbox{length} (\gamma)\geq l \ \wedge\\ \theta \supseteq \tau \ \wedge\ \theta (h(n,\gamma))=0]]\} \end{eqnarray*}  
The set $S$ contains finite strings $\tau$ which say some tree $U_n$ is ill-founded, but for some length $l$, there is no extension 
of $\tau$ in $T$ that says some subtree of $U_n$ with root length at least $l$ is ill-founded.  In short, $\tau$ says $U_n$ is 
ill-founded but there is no sequence of extensions in $T$ to witness it.  

Let $U_n$ be an arbitrary ill-founded tree ($n \notin \KlO$) and let $Z$ be an infinite path through $U_n$.  Then the subtrees extending 
initial segments of $Z$, $U_{h(n,Z|l)}$ for $l \in \omega$, are also ill-founded.  Hence $\KlO(h(n,Z|l)) = 0$ for every $l \in \omega$ 
so $n$ does not witness that $\KlO$ meets $S$.  We conclude that $\KlO$ does not meet $S$\@.  Since S is $\Sigma^0_2 (T)$ and we have 
assumed $\KlO$ is 2-generic relative to $T$, we let $k$ be such that for any $\sigma \in T$ extending $\KlO | k$ we have $\sigma \notin S$\@.  
We can now use the fact that these extensions are sufficiently well behaved to calculate $\KlO$ from $T$.  

\begin{claim} For any number $e$, we have 
$e \in \KlO \iff \neg \exists \sigma \in T[\sigma \supseteq \KlO | k\ \wedge\ \sigma (e) = 0]$. \end{claim}
\begin{proof} $(\Longleftarrow)$ Let $\sigma = \KlO | \max(k,e)+1$.  $\sigma \in T$ since $\KlO$ is a path in $T$, so 
$\sigma (e) = \KlO (e) = 1$.  Hence $e \in \KlO$.

$(\Longrightarrow)$ Let $e \in \KlO$ and suppose the conclusion fails, witnessed by $\sigma$.  $U_e$ is well-founded since 
$e \in \KlO$.  We will construct an infinite path through $U_e$ to get the desired contradiction.  We use an induction to 
simultaneously construct paths $\gamma$ through $U_e$ and $\theta$ through $T$.  Let $j_0$ denote $e$ and $j_{m+1}$ denote 
$h(j_m,\gamma_{m+1})$.  We maintain inductively that $\theta_m (j_m)=0$.  

We begin with $\gamma_0 = \langle \rangle$ and $\theta_0 = \sigma$ and note $\theta_0 (j_0)) = \sigma (e) = 0$ by our assumption.  
Let $\gamma_m$ and $\theta_m$ be given.  $\theta_m \supseteq \sigma \supseteq \KlO |k$ so $\theta_m \notin S$.  Hence we have 
\begin{eqnarray*}  \forall n \forall l[\theta_m (n) \not= 0\ \vee\ \exists \alpha \in U_n \exists \beta \in T[\mbox{length} (\alpha) 
\geq l \ \wedge\ \beta \supseteq \theta_m \ \wedge \\ \beta (h(n,\alpha))=0]] \end{eqnarray*}
Choosing $n=j_m$ and $l=$length$(\gamma_m)+1$ and noting by our induction hypothesis $\theta_m (j_m)=0$, we get 
\begin{eqnarray*} \exists \alpha \in U_{j_m} \exists \beta \in T[\mbox{length} (\alpha) \geq \mbox{length} (\gamma_m) +1 
\ \wedge\ \beta \supseteq \theta_m \ \wedge \\ \beta (h(j_m,\alpha))=0] \end{eqnarray*}
We now let $\gamma_{m+1} = \alpha$ and $\theta_{m+1} = \beta$.  We note that $\gamma_{m+1} \supseteq \gamma_m$ since $\gamma_{m+1} \in 
U_{h(j_{m-1},\gamma_m)}$ and that $\theta_{m+1}(j_{m+1})=\theta_{m+1}(h(j_m,\gamma_{m+1}))=0$, completing the induction.
\end{proof}
Thus $\KlO$ is $\Pi^0_1 (T)$, contradicting $\KlO$ being 2-generic relative to $T$.    
\end{proof}
The next lemma will be used in showing that $\KlO^{(\alpha + 1)}$ is not $n$-generic relative to any perfect tree, given that $\KlO^{(\alpha)}$ is not.
We prove it by applying the same ideas used in the above lemma to the column of $\KlO^{(\alpha)}$ which computes $\KlO$.

\begin{lemma} Let $X \geq_T \KlO$ be 2-generic relative to the perfect tree $T$\@.  Then ${T \geq_T \KlO}$. 
\label{sucO}\end{lemma}
\begin{proof} Let $\Phi$ be a Turing reduction such that $\Phi (X) = \KlO$.  We define $S$ as before, this time for the image under 
$\Phi$.  
\begin{eqnarray*} S = \{ \tau \in T \setsep \exists n \exists l [[\Phi (\tau)](n)=0\ \wedge\ \neg (\exists \gamma \in U_n) (\exists \theta 
\in 
T) [\mbox{length} (\gamma)\geq l\ \wedge\\ \theta \supseteq \tau\ \wedge\ [\Phi (\theta)](h(n,\gamma))=0]]\} \end{eqnarray*} 
We note $S$ is $\Sigma^0_2(T)$ and $X$ does not meet $S$\@.  Since $X$ is 2-generic relative to $T$, we let $k$ be such that for any $\sigma \in 
T$ extending $X|k$ we have $\sigma \notin S$\@.  We now claim that for any $e$, we have $e \in \KlO$ if and only if there does not 
exist a $\sigma \in T$ with $\sigma \supseteq X|k$ and $[\Phi (\sigma)](e)=0$.  This is proved in substantially the same manner as the 
claim in the previous lemma.  As a result, $\KlO$ is $\Pi^0_1 (T)$.  We can now use the fact that generics do not accelerate arithmetic 
definability.  Since $X$ is 2-generic relative to $T$ and $\KlO \leq_T X$, we get $\KlO \leq_T T$ as desired.  
\end{proof} 
  
\begin{corr} For all $n\in\omega$, $\KlO^{(n)}$ is not 2-generic relative to any perfect tree.\label{sublimitO}\end{corr}
\begin{proof} Fix $n$ and suppose not, witnessed by $T$.  We show by induction on $m \leq n$ that $\KlO^{(m)} \leq_T T$\@.  Given 
$\KlO^{(m)} \leq_T T$, we relativize Lemma \ref{sucO} to $\KlO^{(m)}$ to get $\KlO^{(m+1)} \leq_T T$, completing
the induction.  Hence $\KlO^{(n)} \leq_T T$, contradicting our assumption that $\KlO^{(n)}$ is 2-generic relative to $T$.
\end{proof}

\begin{corr} The statement ``All but countably many reals are 2-generic relative to some perfect tree'' fails to hold in $\Pi^1_1$-CA\@. 
\end{corr}
\begin{proof} Consider the standard model of $\Pi^1_1$-CA containing the reals $X$ such that $\exists n [X \leq_T \KlO^{(n)}]$.  The 
set $\{\KlO^{(n)} \setsep n \in \omega \}$ is not a countable set in this model.  
\end{proof}
To handle limit ordinals, we use a lemma in the style of Enderton and Putnam \cite{EndPut}.

\begin{lemma}[Slaman \cite {limitcase}] Let $A$ be a set and $\lambda$ a recursive limit ordinal.  Suppose that for all 
$\beta < \lambda$, $\KlO^{(\beta)} \leq_T A$.  Then $\KlO^{(\lambda)}$ is $\Sigma^0_5 (A)$. \label{limitO}\end{lemma}
\begin{proof} We continue to use $\KlO = \{e \setsep U_e $ is well founded$\}$ where $U_e$ denotes the $e$th recursive tree in 
$\omega^{< \omega}$.  Since $\KlO \leq_T A$ we can define $\KlO$ from $A$ by noting that $U_e$ is well founded iff $U_e$ has no 
infinite path recursive in $A$\@.  Hence $\KlO$ is uniformly $\Pi^0_3 (A)$.  Similarly, we can get $\KlO^\KlO$ is uniformly 
$\Pi^0_4 (A)$ by $X = \KlO \oplus \KlO^\KlO$ iff 
$$(X)_0 = \KlO\ \wedge\ (e \in (X)_1 \leftrightarrow U^{(X)_0}_e \mbox{ has no infinite path recursive in }A)$$ 
where $(X)_0$ and $(X)_1$ denote the two columns of $X$\@.  

We extend this idea to find a uniform definition for $\KlO^{(\lambda)}$.  Fix a system of notations, $o$, for $\lambda$.  We have 
\begin{eqnarray*} & (b,k) \in \KlO^{(\lambda)} \Leftrightarrow \exists m [ \{m\}^A = Y\ \wedge\ k \in (Y)_b\ \wedge\ [\forall c 
\forall d [o(c) < o(b) \rightarrow & \\ & ((o(c)=o(d)+1 \rightarrow \Gamma ((Y)_d, (Y)_c))\ \wedge & \\ & (o(c) 
\mbox{ a limit ordinal } 
\rightarrow \forall n \forall p [(Y)_{c_n}(p) = ((Y)_c)_n (p)]))]]] & \end{eqnarray*}
where $c_0, c_1, c_2, \ldots$ is the fundamental sequence for $o(c)$ and $\Gamma (X,Z)$ is the statement 
$$\forall e[e \in Z \leftrightarrow U^X_e \mbox{ has no infinite path recursive in } A]$$
Then $\Gamma$ is $\Pi^0_4 (A)$ so $\KlO^{(\lambda)}$ is $\Sigma^0_5 (A)$. 
\end{proof}
If we repeat the proof with $\overline{\KlO^{(\lambda)}}$ we improve the result slightly to $\KlO^{(\lambda)}$ is $\Delta^0_5 
(A)$.  

Now we can complete our induction.

\begin{theorem} Let $\lambda$ be the least ordinal such that $\sup_{\beta < \lambda} (\beta$th admissible$) = \lambda$.  Then for all 
$\alpha < \lambda$ we have $\KlO^{(\alpha)}$ is not 5-generic relative to any perfect tree.  \end{theorem}
\begin{proof} Suppose not, witnessed by $\beta$ and $T$.  We define the function $f$ by $f(0)=\omega^{\scriptscriptstyle C \! K}_1$, 
$f(\delta + 1) =$ least admissible greater than $f(\delta)$, and for limit $\delta$, $f(\delta) = \sup_{\xi < \delta} f(\xi)$.  We note 
that $\lambda$ is the least fixed point of $f$.  Using the fact that $\omega^{\KlO^{(\delta)}}_1 < \omega^{\KlO^{(\delta + 1)}}_1$ for 
any $\delta$ \cite{Sacks}, we see by induction that $f(\delta) \leq \omega^{\KlO^{(\delta)}}_1$ for all $\delta$.

Let $\alpha$ be least such that $\KlO^{(\alpha)} \not\leq_T T$.    Then $\alpha \leq \beta < \lambda$ so $\alpha < f(\alpha)$.  If 
$\alpha = \xi + 1$ for some $\xi$ then Lemma \ref{sucO} relativized to $\KlO^{(\xi)}$ would result in a contradiction.  Hence $\alpha$ is 
a nonzero limit ordinal so we can choose $\gamma < \alpha$ such that $\alpha < f(\gamma)$.  Then $\alpha < \omega^{\KlO^{(\gamma)}}_1$ so 
we can fix a system of notations, $o$, for $\alpha$ recursive in $\KlO^{(\gamma + 1)}$.  Since $\gamma + 1 < \alpha$ we have 
$\KlO^{(\gamma + 1)} \leq_T T$\@.  We now apply Lemma \ref{limitO} relativized to $\KlO^{(\gamma + 1)}$ for $T$ to get that 
$\KlO^{(\alpha)} \leq_T T$ for a contradiction.
\end{proof}
\section{1-generics}
In the 1-generic case, we can use a variety of approaches to identify sets of reals that are 1-generic relative to some perfect tree and 
sets whose members cannot have this property.

A real is said to be ranked if it is a member of a countable $\Pi^0_1$ set.  Equivalently, a real is ranked if it is a path through a 
recursive tree with no perfect subtrees.  The reader is referred to Cenzer et al.\ \cite{Ranked} for details on the topic, including a 
proof that for all recursive ordinals $\alpha$ there is a ranked set of degree $\mathbf{0^{(\alpha)}}$.  Here we demonstrate these 
reals are not 1-generic relative to any perfect tree.

\begin{prop} If $X$ is 1-generic relative to some perfect tree, then $X$ is not ranked. \end{prop}
\begin{proof} Suppose not.  Let $X$ be 1-generic relative to the perfect tree $T$ and a path through the recursive tree $U$ with no 
perfect subtrees.  Let $S = \{ \sigma \in T \setsep \sigma \notin U \} $.  Then $S$ is recursive in $T$ and $X$ does not meet $S$, so there 
exists an $n$ such that no $\tau \in T$ extending $X|n$ is in $S$\@.  Hence for every $\tau \in T$ such that $\tau \supseteq X|n$ we 
have $\tau \in U$\@.  But then $U$ has a perfect subtree, for a contradiction.  
\end{proof}
It follows from Cenzer et al.\ \cite{Ranked} and this lemma that there are reals arbitrarily high in the hyper\-arithmetic degrees which 
are not 1-generic relative to any perfect tree.  We note the proof of Lemma \ref{mainlem} can be relativized to start with a perfect tree in 
place of $2^{<\omega}$.  Using this we observe that no Turing degree can contain both a ranked set and a real 2-generic relative to some 
perfect tree (and no truth table degree a ranked set and a real 1-generic relative to some perfect tree).  Hence no $\Delta^0_2$ set is 
2-generic relative to some perfect tree, and the degrees $\mathbf{0^{(\alpha)}}$ for any recursive $\alpha$ contain no reals 2-generic 
relative to some perfect tree.  

We can also attempt to classify which reals are 1-generic relative to some perfect tree by use of the r.e. (Ershov) and REA hierarchies.  
The reader is referred to Jockusch and Shore \cite{REA} for details on these hierarchies.  We begin by observing that no real whose 
degree is at a finite level of the REA hierarchy (hence also the r.e. hierarchy) is 1-generic relative to some perfect tree.

\begin{prop}[Slaman \cite{nREA}] Let $n\in\omega$, $X$ a real of $n$-REA degree.  Then $X$ is not 1-generic 
relative to any perfect tree. \end{prop}
\begin{proof} Fix $n$ and $X$ and let $W$ be an $n$-REA set with $X \equiv_T W$\@.  Let $W_1, W_2, \ldots \\ W_n = W$ witness that
$W$ is $n$-REA; for all $i \leq n$ we have $W_i \leq_T W_{i+1}$ and $W_{i+1}$ is r.e.($W_i$).  Suppose $X$ is 1-generic relative to $T$.  We show by induction that for all $m \leq n$ we have $W_m \leq_T T$ 
using the following claim:   

\begin{claim} Let $Y$ be r.e.($T$) and $Y \leq_T X$\@.  Then $Y \leq_T T$. \end{claim}
\begin{proof} It suffices to show $\overline{Y}$ is r.e.($T$)\@.  Since $Y \leq_T X$, let $\{e\}^X = Y$\@.  Let $$S=\{q \setsep \exists n [ 
\{e\}^q (n)\!\downarrow = 0 \wedge n\in Y ] \}$$  We note $X \notin S$ and $S$ is r.e.($T$) since Y is r.e.($T$)\@.  Hence for some 
$l$, every $q$ extending $X|l$ is not in $S$\@.  We can now describe $\overline{Y}$ by noting that $n \in \overline{Y}$ iff 
\\$\exists q 
\supseteq X|l [\{e\}^q (n)\!\downarrow = 0]$.  Hence $\overline{Y}$ is r.e.($T$) as desired.  
\end{proof}
For the induction, given $W_m \leq_T T$ we note that $W_{m+1}$ is r.e.($W_m$), hence r.e.($T$), and apply the claim to $W_{m+1}$\@.  As 
a result $W \leq_T T$ so $X \leq_T T$ for a contradiction.  
\end{proof}
We might next hope to show sets of $\omega$-REA degree are not 1-generic relative to any perfect tree.  However, we cannot even
do this for sets which are $\omega$-r.e.  In proving the Friedberg Inversion Theorem for the truth table degrees, J. Mohrherr 
\cite{ttdegrees} showed by a reduction that there is a 1-generic $G$ such that $G \leq_{tt} 0^\prime$, hence $G$ is $\omega$-r.e.  
Here we provide a direct construction.  We use the definition \cite{REA} that $X$ is $\omega$-r.e.\ if for some partial recursive 
$\psi : \omega \times \omega \to 2$ we have $X(n) = \psi (b,n)$ where b is least such that $\psi (b,n)\!\downarrow$.

\begin{prop} There is a 1-generic real which is $\omega$-r.e. \end{prop}
\begin{proof} In this construction we extend to meet the first r.e. set we find while still looking for earlier r.e. sets skipped over.  
If we find a set that has been skipped, we start over again from that point.  To limit the number of injuries, we require the $n$th 
r.e.\ set to extend the first $n$ bits of the current string.

We build in stages our generic $X \in 2^\omega $ and the partial recursive witness, $\psi$, that $X$ is 
$\omega$-r.e.  We also use some numeric variables for bookkeeping.  $r$ denotes the r.e.\ set we are looking at, $k_i$ for $i \in 
\omega$ the number of corrections at the $i$-th r.e.\ set, and $m_i$ for $i \in \omega$ the length of the initial segment of $X$ 
currently meeting the $i$-th r.e.\ set (0 if not yet met).  We start with $\psi = \emptyset$, $X=\langle \rangle$, $r=0$, and 
$k_i,m_i = 0$ for all $i$.  

At stage $n+1$ we search simultaneously for $i$ such that $m_i = 0$, $\tau \supseteq X_n|\max(i,\max_{j<i} (m_j))$ with 
length$(\tau) > n$, and $s$ to find $\{i\}^\tau_s \!\downarrow$.  If $i>r$ we have found a new r.e.\ set and add it by letting 
$m_i = $length$(\tau)$, $X_{n+1} = \tau$, $r=i$, and for $l$ such that length$(X_n) < l \leq $length$(X_{n+1})$ we let 
$\psi (2^l + l -k_l,l) = X_{n+1} (l)$.  If instead $i<r$ we have found a r.e.\ set we have skipped over and restart at that point.  We do 
this by first setting $m_i = $length$(\tau)$ and for $l$ such that $\max(i,\max_{j<i} (m_j)) < l \leq m_r$ setting $k_l$ to $k_l + 1$.  
We next reset $X_{n+1}$ to $\tau$ (it will not extend $X_n$, but will extend $X_n| \max(i,\max_{j<i}(m_j))$).  Then, for $l$ such that 
$\max(i,\max_{j<i} (m_j)) < l \leq m_i$ we extend $\psi$ by letting $\psi (2^l + l - k_l,l) = X_{n+1}(l)$.  Finally, for $j$ with 
$i < j \leq r$ we set $m_j = 0$ and last set $r=i$.

Since we require the $n$th r.e.\ set to extend $X|n$, after stage $i$, the value of $X(i)$ can only be changed when $X$ is changed to 
meet the $n$th r.e.\ set for some $n<i$.  By the usual Friedberg-Muchnik counting of injuries, $X(i)$ can be changed at most $2^i$ 
times after stage $i$, so at most $2^i + i$ times in all.  Hence $\psi(b,n)$ witnesses $X$ is $\omega$-r.e.\ by 
starting with 
$b = 2^i+i$ and moving $b$ down one every time a correction is made.  We see that the resulting $X$ is 1-generic since if 
$X$ does 
not meet the $i$th r.e.\ set then there is no extension of $X|\max(i,\max_{j<i}(m_j))$ which meets the $i$th r.e.\ set.
\end{proof}
We note that by the REA Completeness Theorem (Jockusch and Shore \cite{REA}) this gives that for every $X \geq_T \emptyset^{(\omega)}$ 
there are sets $A$ and $J$ such that $X \equiv_T A \oplus J$ where $J$ is 1-generic ($A$).  

There is still considerable room left to explore in determining which reals are 1-generic relative to some perfect tree.  In particular,
it is not yet known if every real not 1-generic relative to any perfect tree is hyper\-arithmetic.
\nocite{*}
\bibliography{GenericRPT}
\end{document}